\documentclass[12pt]{article}
\usepackage{amsfonts}
\usepackage{bbm}
\usepackage{mathrsfs}
\usepackage{amssymb}

\title{Recognizing The Semiprimitivity of $\mathbb{N}$-graded \\
Algebras via Gr\"obner Bases\thanks{Project supported by the
National Natural Science Foundation of China (10971044).}}

\author{Huishi Li\\
{\small Department of Applied Mathematics}\\
{\small College of Information Science and Technology}\\
{\small Hainan University}\\
{\small  Haikou 570228, China}}

\date{}

\begin{document}
\maketitle
\begin{center}
\begin{minipage}{120mm}
{\small {\bf Abstract.}  Let $K\langle X\rangle =K\langle X_1,\ldots
,X_n\rangle$ be the free $K$-algebra on $X=\{ X_1,\ldots ,X_n\}$
over a field $K$, which is equipped with a weight
$\mathbb{N}$-gradation (i.e., each $X_i$ is assigned a positive
degree), and let ${\cal G}$ be a finite homogeneous Gr\"obner basis
for the ideal $I=\langle{\cal G}\rangle$ of $K\langle X\rangle$ with
respect to some monomial ordering $\prec$ on $K\langle X\rangle$. It
is proved that if the monomial algebra $K\langle X\rangle /\langle
{\bf LM}({\cal G})\rangle$ is semi-prime, where ${\bf LM}({\cal G})$
is the set of leading monomials of ${\cal G}$ with respect to
$\prec$, then the $\mathbb{N}$-graded algebra $A=K\langle X\rangle
/I$ is semiprimitive (in the sense of Jacobson). In the case that
${\cal G}$ is a finite non-homogeneous Gr\"obner basis with respect 
to a graded monomial ordering $\prec_{gr}$, and the 
$\mathbb{N}$-filtration $FA$ of the algebra $A=K\langle X\rangle /I$ 
induced by the $\mathbb{N}$-grading filtration $FK\langle X\rangle$ 
of $K\langle X\rangle$ is considered, if the monomial algebra 
$K\langle X\rangle /\langle {\bf LM}({\cal G})\rangle$ is 
semi-prime, then it is proved that the associated 
$\mathbb{N}$-graded algebra $G(A)$ and the Rees algebra 
$\widetilde{A}$ of $A$ determined by $FA$ are all semiprimitive. }
\end{minipage}\end{center} {\parindent=0pt\par

}

\renewcommand{\thefootnote}{\fnsymbol{footnote}}
\setcounter{footnote}{-1}


\def\NZ{\mathbb{N}}
\def\QED{\hfill{$\Box$}}
\def \r{\rightarrow}
\def\mapright#1#2{\smash{\mathop{\longrightarrow}\limits^{#1}_{#2}}}

\def\v5{\vskip .5truecm}
\def\OV#1{\overline {#1}}
\def\hang{\hangindent\parindent}
\def\textindent#1{\indent\llap{#1\enspace}\ignorespaces}
\def\item{\par\hang\textindent}

\def\LH{{\bf LH}}\def\LM{{\bf LM}}\def\LT{{\bf
LT}}\def\KX{K\langle X\rangle} \def\KS{K\langle X\rangle}
\def\B{{\cal B}} \def\LC{{\bf LC}} \def\G{{\cal G}} \def\FRAC#1#2{\displaystyle{\frac{#1}{#2}}}
\def\SUM^#1_#2{\displaystyle{\sum^{#1}_{#2}}} \def\T#1{\widetilde #1} \def\HL{{\rm LH}}

\section*{1. Introduction}
Let $K\langle X\rangle =K\langle X_1,\ldots ,X_n\rangle$ be the free
$K$-algebra on $X=\{ X_1,\ldots ,X_n\}$ over a field $K$, which is
equipped with a weight $\mathbb{N}$-gradation (i.e., each $X_i$ is
assigned a positive degree), and let ${\cal G}$ be a finite
Gr\"obner basis for the ideal $I=\langle{\cal G}\rangle$ of
$K\langle X\rangle$ with respect to some monomial ordering $\prec$
on $K\langle X\rangle$. Consider the algebra $A=\KX /I$, the
associated monomial algebra $\OV A=\KX /\langle\LM (\G )\rangle$ of
$A$ (where $\LM (\G )$ is the set of leading monomials of $\G$ with
respect to $\prec$), the $\mathbb{N}$-filtration $FA$ of $A$ induced 
by the $\mathbb{N}$-grading filtration of $K\langle X\rangle$, the 
associated $\NZ$-graded algebra $G(A)$ and the Rees algebra $\T A$ 
of $A$ determined by $FA$ (see Section 3 for the definitions). In 
[Li2] and [Li3], it has been proved that many structural properties 
of $\OV A$ can be transferred to $A$, $G(A)$ and $\T A$ (see [Li4] 
for more details). In this paper, we first show that if $\G$ is a 
finite homogeneous Gr\"obner basis and if the monomial algebra $\OV 
A$ is semi-prime,  then the $\mathbb{N}$-graded algebra $A=K\langle 
X\rangle /I$ is semiprimitive (in the sense of Jacobson). In the 
case that ${\cal G}$ is a finite non-homogeneous Gr\"obner basis 
with respect to a  graded monomial ordering $\prec_{gr}$, if $\OV A$ 
is semi-prime, then we show that the $\mathbb{N}$-graded algebras 
$G(A)$ and $\widetilde{A}$ are all semiprimitive. Since the 
semi-primeness of the monomial algebra $\OV A$ can be determined in 
an algorithmic way ([G-IL], [G-I]), our results are algorithmically 
realizable in case the algorithms given in loc. cit. are implemented 
on computer. \v5

Throughout this paper, $K$ denotes a field, algebras considered are
associative $K$-algebras with multiplicative identity 1, and ideals
considered in an algebra are meant two-sided ideals. For a subset
$U$ of an algebra $A$, we write $\langle U\rangle$  for the ideal
generated by $U$ in $A$. Moreover, we use $\NZ$, respectively
$\mathbb{Z}$, to denote the set of nonnegative integers,
respectively the set of integers.

\section*{2. Some Known Results on Monomial Algebras}
Let $\KX =K\langle X_1,\ldots ,X_n\rangle$ be the free $K$-algebra
on $X=\{ X_1,\ldots ,X_n\}$, and $\B =\{
X_{i_1}^{\alpha_1}X_{i_2}^{\alpha_2}\cdots
X_{i_s}^{\alpha_s}~|~X_{i_j}\in X,~\alpha_j\in\NZ\}$ the standard
$K$-basis of $\KX$ consisting of all monomials (words) in
$X_{i_j}$'s. For convenience, we use lowercase letters $w, u, v, s,
\ldots$ to denote monomials in $\B$. In this section we recall from
[G-IL] and [G-I] how to recognize the Jacobson semiprimitivity of a
finitely presented monomial algebra $R=\KX /\langle\Omega\rangle$ in
a computational way, where $\Omega =\{ u_1,\ldots ,u_s\}\subset\B$
is a reduced finite subset of monomials (see the definition below)
such that $\Omega\cap X=\emptyset$. \v5

For $u,v\in\B$, we say that $v$ {\it divides} $u$, denoted by $v|u$,
if $u=wvs$ for some $w,s\in\B$.  We say that a subset
$\Omega\subset\B$ is {\it reduced} if $v,u\in\Omega$ and $v\ne u$
implies $v{\not |}u$.

{\parindent=0pt\v5

{\bf 2.1. Theorem} Let $\Omega =\{ u_1,\ldots ,u_s\}$ be a reduced
finite subset of $\B-X$ and $R=\KX /\langle\Omega\rangle$. \par

(i) ([G-IL], Theorem 16) The Jacobson radical $J(R)$ of $R$
coincides with the upper nilradical Nil$(R)$ of $R$.\par

(ii) ([G-I], Theorem 2.27) $R$ is semiprimitive (in the sense of
Jacobson), i.e., $J(R)=\{ 0\}$, if and only if $R$ is semi-prime,
i.e., $a\in R$ and $aRa=\{ 0\}$ implies $a=0$.\par\QED} \v5

Let $\Omega =\{ u_1, \ldots ,u_s\}$ be a  reduced finite subset of
$\B$, and $\langle\Omega\rangle$ the monomial ideal generated by
$\Omega$. Then the set of {\it normal monomials} (mod
$\langle\Omega\rangle$) in $\B$ is defined as
$$N(\Omega )=\{ w\in\B~|~u{\not |}w,~u\in\Omega\} .$$
For each $u_i\in\Omega$, say $u_i=X_{i_1}^{\alpha_1}\cdots
X_{i_r}^{\alpha_r}$ with $X_{i_j}\in X$ and $\alpha_j\in\NZ$, we
write $l(u_i)=\alpha_1+\cdots +\alpha_r$ for the length of $u_i$ .
Put
$$\ell =\max\left\{\left. l(u_i)~\right |~u_i\in\Omega\right\} .$$
Then the {\it Ufnarovski graph} of $\Omega$ (in the sense of [Uf1]),
denoted by $\Gamma (\Omega)$, is defined as a directed graph, in
which the set of vertices $V$ is given by
$$V=\left\{\left. v_i\in N(\Omega )~\right |~l(v_i)=\ell -1\right\} ,$$
and the set of edges $E$ contains the edge $v_i\rightarrow v_j$ if
and only if there exist $X_k$, $X_t\in X$ such that 
$v_iX_k=X_tv_j\in N(\Omega )$.  Since $\Omega$ is finite, the 
directed graph $\Gamma (\Omega )$ is thereby practically 
constructible.\par

{\parindent=0pt\v5

{\bf Remark} (i) Note that we have defined the number $\ell$ above
as in [G-IL]. While in [G-I] this number was defined as
$m+1=\max\left\{ l(u_i)~\Big |~u_i\in\Omega\right\}$. So, in the
subsequent results we shall use $\ell$ and $\ell -1$ instead of
$m+1$ and $m$.\par

(ii) To better understand the practical application of $\Gamma
(\Omega )$, it is essential to notice that a Ufnarovski graph is
defined by using the {\it length $l(u)$ of the monomial} ({\it
word}) $u\in \B$ instead of using the {\it degree of $u$ as a 
homogeneous element in $\KS$} whenever a weight $\NZ$-gradation of 
$\KS$ is used (see Section 3), though both notions coincide when 
each $X_i$ is assigned the degree 1.}\v5

Basic notions from classical graph theory fully suit an Ufnarofski
graph $\Gamma (\Omega )$. For instance, a {\it route of length} $m$
 in $\Gamma (\Omega )$ with $m\ge 1$ is a sequence of edges
$$v_0\rightarrow v_1\rightarrow v_2\rightarrow\cdots\rightarrow v_{m-1}\rightarrow v_m.$$
If in a route no edge appears repeatedly then it is called a {\it
simple route}. A simple route with $m\ge 1$ and $v_0=v_m$ is called
a {\it cycle}. If, as an undirected graph, there is a route between
any two distinct vertices of $\Gamma (\Omega )$, then $\Gamma
(\Omega )$ is called {\it connected}.  A {\it connected component}
of $\Gamma (\Omega )$ ( as an undirected graph) is a connected
subgraph which is connected to no additional vertices. \par

Let $\Gamma (\Omega )$ be as above. It follows from [Uf1] that there
is a one-to-one correspondence between the monomials (words) of
length $\ge\ell -1$ in $N(\Omega )$ and the routes in $\Gamma
(\Omega )$, that is, if $u\in N(\Omega )$ and $u=X_{i_1}\cdots
X_{i_t}$ with $t\ge \ell -1$, then $u$ is mapped to the route
$$\mathscr{R}(u):\quad v_0\rightarrow v_1\rightarrow\cdots\rightarrow v_m,$$
where\par
$$\begin{array}{l} m=t-\ell +1,~\hbox{and}\\
v_j=X_{i_{j+1}}X_{i_{j+2}}\cdots X_{i_{j+\ell -1}},~0\le j\le
m.\end{array}$$ In [G-I] the following notions are introduced. A
vertex of $\Gamma (\Omega )$ is called a {\it cyclic vertex} if it
belongs to a cycle. A normal monomial $u\ne 1$ (i.e., $u\in N(\Omega
)-\{ 1\}$) is called a {\it cyclic normal monomial} if $l(u)\le \ell
-1$ and $u$ is a suffix of some cyclic vertx of $\Gamma (\Omega )$,
or, if $l(u)>\ell -1$ and the associated route $\mathscr{R}(u)$ of
$u$ is a subroute of some cyclic route.{\parindent=0pt\v5

{\bf 2.2. Proposition}  Let $\Omega =\{ u_1,\ldots
,u_s\}\subset\B-X$ be reduced, $R=\KX /\langle\Omega\rangle$, and
let $J(R)$ be the Jacobson radical of $R$. If $u\in N(\Omega )-\{
1\}$, then the following statements hold.
\par

(i) ([G-I], Lemma 2.13) $u$ is cyclic if and only if there exists a
monomial $v\in\B$ such that $(uv)^q\not\in\langle\Omega\rangle$ for
all integer $q\ge 0$.\par

(ii) ([G-I], Corollary 2.17) $\OV u\in J(R)$ if and only if $u$ is
noncyclic, where $\OV u$ is the residue class representaed by $u$ in
$R$.\v5

{\bf 2.3. Theorem} ([G-I], Theorem 2.21) Let the monomial algebra
$R=\KX /\langle\Omega\rangle$ be as in Proposition 2.2. Then $R$ is
semiprimitive if and only if any monomial $u\in N(\Omega )$ with
$1\le l(u)\le \ell$ is cyclic.\par \QED\v5

{\bf Remark} (i) By Theorem 2.1(ii), it is clear that if $R=\KX
/\langle\Omega\rangle$ is a prime ring, then $R$ is
semiprimitive.\par

(ii) The reader is referred to [G-IL] and [G-I] for the algorithms
written for determining the semi-primeness and the primeness (and
hence the semiprimitivity) of a finitely presented monomial algebra
$R=\KX /\langle\Omega\rangle$.}\v5

\section*{3. The Main Results}
In this section, we prove the main results of this paper (Theorem
3.2, Theorem 3.3, Theorem 3.5).  The Gr\"obner basis theory for 
ideals in a free $K$-algebra is referred  to ([Ber], [Mor], [Gr], 
[Uf2]). \v5

To begin with, let $K$ be a field and $\KX =K\langle X_1,\ldots
,X_n\rangle$ the free $K$-algebra on $X=\{ X_1,\ldots ,X_n\}$. As
before the standard $K$-basis of $\KX$ is denoted by $\B$. We fix a
weight $\NZ$-gradation for $\KX$, that is, $\KX
=\oplus_{p\in\NZ}\KX_p$ in which, each $X_i$ has an assigned
positive degree $m_i$, $1\le i\le n$, and the degree-$p$ homogeneous
part $\KX_p$ is the $K$-vector space spanned by all monomials of
degree $p$. For a nonzero homogeneous element $H\in \KX_p$, we write 
$d(H)$ for the degree of $H$, i.e., $d(H)=p$. Note that every 
monomial $w\in\B$ is a homogeneous element. If $I$ is a graded ideal 
of $\KX$ (i.e., $I$ is generated by homogeneous elements), then 
$A=\KX /I$ is an $\NZ$-graded algebra, that is, 
$A=\oplus_{p\in\NZ}A_p$ with the degree-$p$ homogeneous part 
$A_p=(\KX_p+I)/I$.\par

Moreover, let $\prec$ be a monomial ordering on $\B$, i.e., $\prec$
is a well-ordering on $\B$ such that $u\prec v$ implies $wus\prec
wvs$,  and $v=wus$ with $w\ne 1$ or $s\ne 1$ implies $u\prec v$, for
all $w,u,v,s\in\B$. With the monomial ordering $\prec$ fixed on 
$\B$, each subset $S$ of $\KX$ is associated to a subset of 
monomials $\LM (S)=\{ \LM (f)~|~f\in S\}\subset\B$, where if $f\in 
S$ and $f=\sum_{i=1}^s\lambda_iw_i$ with $\lambda_i\in K$ and 
$w_i\in\B$ such that $w_1\prec w_2\prec\cdots\prec w_s$, then $\LM 
(f)=w_s$. $\LM (S)$ is usually referred to as the {\it set of 
leading monomials} of $S$. By the classical Gr\"obner basis theory 
of $\KX$, in principle every nonzero ideal $I$ of $\KX$ has a 
nontrivial (finite or infinite) Gr\"obner basis $\G$ in the sense 
that $\G$ is a proper subset of $I-\{ 0\}$ and $\langle\LM 
(I)\rangle =\langle\LM (\G )\rangle$. A Gr\"obner basis $\G$ 
consisting of homogeneous elements of $\KX$ is called a {\it 
homogeneous Gr\"obner basis}. \v5

In proving our main theorems, we need a fundamental result
concerning the Jacobson radical of a $\mathbb{Z}$-graded ring, which
is due to G. Bergman (cf. [Row]). {\parindent=0pt\v5

{\bf 3.1. Theorem} Let $R=\oplus_{n\in\mathbb{Z}}R_n$ be a
$\mathbb{Z}$-graded ring and $J(R)$ the Jacobson radical of $R$.
Then\par

(i) $J(R)$ is a graded ideal of $R$; and\par

(ii) if $n\ne 0$ and $a\in R_n$, then $1+a$ is invertible if and
only if $a$ is nilpotent.\par\QED\v5

{\bf 3.2. Theorem} Let $\KX$ and the monomial ordering $\prec$ on
$\B$ be as fixed above, and let $\G =\{ g_1,\ldots ,g_s\}$ be a
finite homogeneous Gr\"obner basis for the ideal
$I=\langle\G\rangle$, such that $\LM (\G )\cap X=\emptyset$ and $\LM
(\G )$ is reduced (in the sense of Section 2). If the monomial
algebra $\OV A=\KX /\langle\LM (\G )\rangle$ is semi-prime, then the
$\NZ$-graded algebra $A=\KX /I$ is semiprimitive. \vskip 6pt

{\bf Proof} Let $N(I)$ be the set of normal monomials (mod $I$) in
$\B$, i.e., $N(I)=\{ w\in\B~|~w\not\in\langle\LM (I)\rangle\}$.
Then, since $\G$ is a Gr\"obner basis of $I$, $N(I)=\{ u\in\B~|~\LM
(g_i){\not |u},~g_i\in\G\}$. If, with respect to the $\NZ$-gradation
of $\KX$, $H\in\KX_p$ is a homogeneous element of degree $p$, and if
$H\not\in I$, then by the division by the homogeneous Gr\"obner
basis $\G$, $H$ has a representation
$H=\sum_{i,j}\lambda_{ij}w_{ij}g_jv_{ij}+\sum_t\mu_tu_t$, where
$\lambda_{ij},\mu_t\in K-\{ 0\}$, $w_{ij},v_{ij}\in\B$, $g_j\in\G$,
and $u_t\in N(I)$ such that $d(H)=p=d(u_t)$ for all $t$. Putting
$H'=\sum_t\mu_tu_t$ and considering the nonzero homogeneous element
$\OV H$ of degree $p$ represented by $H$ in $A_p=(\KX_p +I)/I$, we
have
$$\OV H=\OV{H'}=\sum_t\mu_t\OV{u_t},\quad
\mu_t\in K-\{ 0\},~u_t\in
N(I)~\hbox{with}~d(u_t)=p=d(H).\eqno{(1)}$$}\par

Let $J(A)$ be the Jacobson radical of  the $\NZ$-graded algebra
$A=\KX /I$. If $J(A)\ne \{ 0\}$, then it follows from Theorem 3.1
that $J(A)$ is a graded ideal of $A$.  Taking a nonzero homogeneous
element of $J(A)$, say $\OV H\in J(A)\cap A_p$, where $A_p=(\KX_p
+I)/I$ and $H\in\KX_p$, we may replace $\OV H$ by $\OV{H'}$ as in
(1) above. Without loss of generality we assume that $\LM (H')=u_1$
with respect to $\prec$. Our aim below is to show that
{\parindent=.7truecm\vskip 6pt

\item{$(\bullet)$} the normal ~monomial $u_1$ is noncyclic, thereby~
$\OV{u_1}\in J(\OV A)$ by Proposition 2.2, ~where\par $\OV{u_1}$ is
the residue class represented by $u_1$ in $\OV A$. }\vskip 6pt

Assume the contrary that $u_1$ is cyclic (see Section 2 for the
definition). Then, by Proposition 2.2, there is a monomial $v\in\B$
such that $(u_1v)^q\not\in\langle\LM (I)\rangle =\langle\LM (\G
)\rangle$, or equivalently, $(u_1v)^q\in N(I)$ for all  $q\in\NZ$. 
Since $\LM (H'v)=u_1v$, it turns out that 
$$\LM\left ((H'v)^q\right )=\left (\LM (H'v)\right )^q=(u_1v)^q\in N(I),~q\in\NZ.\eqno{(2)}$$
On the other hand, writing $\OV v$ for the residue class represented
by $v$ in $A$, we have $\OV{H'v}=\OV{H'}\OV v\in J(A)$. As
$\OV{H'v}$ is again a homogeneous element of $J(A)$, it follows from
Theorem 3.1 that $\OV{H'v}^m=0~\hbox{for some integer}~m>0$. Hence, 
$(H'v)^m\in I$ and this gives rise to 
$$(u_1v)^m=\left (\LM (H'v)\right )^m=\LM\left ((H'v)^m\right )\in\langle\LM (I)\rangle=
\langle\LM (\G )\rangle.\eqno{(3)}$$ Clearly, $(3)$ contradicts 
$(2)$. Therefore, $u_1$ is noncyclic and consequently $\OV{u_1}\in 
J(\OV A)$, proving the claim $(\bullet)$. \par

Finally, suppose that the monomial algebra $\OV A=\KX /\langle\LM
(\G )\rangle$ is semi-prime. Then it follows from Theorem 2.1(ii)
that  $\OV A$ is semiprimitive and hence $J(\OV {A})=\{ 0\}$. By the
above argument, we conclude that $J(A)=\{ 0\}$; otherwise, by the
claim $(\bullet )$, there would be a normal monomial $u_1\in N(I)$ 
such that $\OV{u_1}\in J(\OV A)=\{ 0\}$ and hence, $u_1\in 
\langle\LM (\G )\rangle$, which is a contradiction. This shows that 
$A$ is semiprimitive. \QED\v5

Let $I$ be an {\it arbitrary} proper ideal of $\KX$, $A=\KX /I$.
Consider the natural $\NZ$-grading filtration $F\KX=\{
F_p\KX\}_{p\in\NZ}$ of $\KX$ determined by a fixed weight
$\NZ$-gradation for $\KX$, that is, $F_p\KX =\oplus_{q\le p}\KX_q$ 
for $p\in\NZ$.  Then $A$ has the natural $\NZ$-filtration $FA=\{ 
F_pA\}_{p\in\NZ}$ induced by the $\NZ$-grading filtration $F\KX$, 
where $F_pA=(F_p\KX+I)/I$ for $p\in\NZ$, and  $A$ has the associated 
$\NZ$-graded algebra $G(A)=\oplus_{p\in\NZ}G(A)_p$ with 
$G(A)_p=F_pA/F_{< p}A$, where $F_{< p}A=\cup_{q<p}F_{q}A$ 
(conventionally we put $F_{< 0}A=\{ 0\})$. In case $I$ is a graded 
ideal of $\KX$,  it is clear that $G(A)\cong A$ as $\NZ$-graded 
algebras. \par

We also recall that any total ordering $\prec$ on $\B$ induces an
ordering $\prec_{gr}$ on $\B$ subject to the rule: For $u,v\in\B$,
$$u\prec_{gr}v\Leftrightarrow d(u)<d(v)~\hbox{or}~ d(u)=d(v)~\hbox{and}~u\prec v,$$
where $d(~)$ is the degree function on the homogeneous elements of 
$\KX$. If $\prec_{gr}$ is a monomial ordering on $\B$, then 
$\prec_{gr}$ is called an {\it $\NZ$-graded monomial ordering}, for 
instance, the commonly used $\NZ$-graded lexicographic ordering on 
$\B$. {\parindent=0pt\v5

{\bf 3.3. Theorem} With notation as before, let $\prec_{gr}$ be an
$\NZ$-graded monomial ordering on $\B$ with respect to a fixed
weight $\NZ$-gradation of $\KX$, and let $\G =\{ g_1,\ldots ,g_s\}$
be a finite (but not necessarily homogeneous) Gr\"obner basis for
the ideal $I=\langle\G\rangle$, such that $\LM (\G )\cap
X=\emptyset$ and $\LM (\G )$ is reduced (in the sense of Section 2).
Consider the algebra $A=\KX /I$.  If the monomial algebra $\OV A=\KX
/\langle\LM (\G )\rangle$ is semi-prime, then $G(A)$ is
semiprimitive. \vskip 6pt

{\bf Proof} Let $\LH (\G )=\{ \LH (g_i)~|~g_i\in\G\}$ be the set of
$\NZ$-leading homogeneous elements of $\G$ with respect to the fixed
weight $\NZ$-gradation of $\KX$, that is, $\LH (g_i)=H_p$ if
$g_i=H_0+H_1+\cdots +H_p$ with $H_j\in\KX_j$ and $H_p\ne 0$. Since
we are using an $\NZ$-graded monomial ordering $\prec_{gr}$, it
follows from ([LWZ], Theorem 2.3.2) or ([Li2], Proposition 3.2) that
$\LH (\G )$ is a Gr\"obner basis for the graded ideal $\langle\LH
(\G )\rangle$ of $\KX$, and  that
$$G(A)\cong\KX /\langle\LH (\G )\rangle$$ as $\NZ$-graded
algebras. Furthermore, under the $\NZ$-graded monomial ordering
$\prec_{gr}$ we have $\LM (\G )=\LM (\LH (\G ))$. Hence  the
$\NZ$-graded algebra $\KX /\langle\LH (\G )\rangle$ has the
associated monomial algebra $\OV A=\KX /\langle\LM (\G )\rangle$.
Consequently, our assertion follows from Theorem 3.2 and the
isomorphism given above.\QED}\v5

Now, we turn to the Rees algebra $\T A$ of the $\NZ$-filtered
algebra $A=\KX /I$, where $I$ is an {\it arbitrary} proper ideal of
$\KX$, the $\NZ$-filtration $FA=\{ F_pA\}_{p\in\NZ}$ for $A$ is as 
constructed before Corollary 3.3, and  $\T A$ is defined as the 
$\NZ$-graded algebra $\T A =\oplus_{p\in\NZ}F_pA$ with the 
multiplication induced by $F_pAF_qA\subseteq F_{p+q}A$ for all 
$p,q\in\NZ$. The relations between $A$, $G(A)$ and $\T A$ are given 
by the algebra isomorphisms $A\cong \T A/\langle 1-Z\rangle$ and 
$G(A)\cong\T A/\langle Z\rangle$, where $Z$ is the homogeneous 
element of degree $1$ in $\T A_1=F_1A$ represented by the 
multiplicative identity element 1 of $A$. Because of these 
relations, the structure of $\T A$ is closely related to the study 
of the homogenized algebra of an algebra defined by relations, the 
regular central extension  and the PBW-deformation of an 
$\NZ$-graded algebra defined by relations (cf. [LWZ], [Li1], [Li2]).
\par\def\KXT{K\langle X,T\rangle}

Consider the free $K$-algebra $\KXT =K\langle X_1,\ldots
,X_n,T\rangle$ in which each $X_i$ has the same positive degree as
fixed in $\KX$, and we assign $d(T)=1$. Write $\T{\B}$ for the
standard $K$-basis for $\KXT$. If $\prec_{gr}$ is some $\NZ$-graded
lexicographic ordering on the standard $K$-basis $\B$ of $\KX$, then
$\prec_{gr}$ extends to a $\NZ$-graded lexicographic ordering
$\prec_{_{T\hbox{-}gr}}$ on $\T{\B}$ subject to
$T\prec_{_{T\hbox{-}gr}} X_i$, $1\le i\le n$. If $f\in\KX$ has the
linear representation $f=\lambda\LM (f)+\sum_i\lambda_iw_i$ with 
$\lambda$, $\lambda_i\in K-\{ 0\}$, $\LM (f)\in\B\cap\KX _p$, 
$w_i\in\B\cap\KX_{q_i}$, then the {\it non-central homogenization} 
of $f$ with respect to $T$ is the homogeneous element
$$\T f=\LC (f)\LM (f)+\sum_i\lambda_iT^{p-q_i}w_i\in\KXT_p .$$
Clearly, $\LM (f)=\LM (\T f)$, where $\LM (f)$ is taken with respect
to $\prec_{gr}$ on $\B$ and $\LM (\T f)$ is taken  with respect to
$\prec_{_{T\hbox{-}gr}}$ on $\T{\B}$. If $I$ is an ideal of $\KX$, 
then we put
$$\T I=\{ \T f~|~f\in I\}\cup\{ X_iT-TX_i~|~1\le i\le
n\},$$ and call $\langle\T I\rangle$, the graded ideal of $\KXT$ 
generated by $\T I$, the {\it non-central homogenization ideal} of 
$I$ in $\KXT$ with respect to $T$.{\parindent=0pt\v5

{\bf 3.4. Proposition} With notation as fixed above, let $I$ be an
ideal of $\KX$ and $\G\subset I$. The following statements are
equivalent.\par

(i) $\G$ is a Gr\"obner basis for $I$ with respect to $\prec_{gr}$
on $\B$.\par

(ii) $\T{\G}=\{\T g~|~g\in\G\}\cup\{ X_iT-TX_i~|~1\le i\le n\}$ is a
homogeneous Gr\"obner basis for $\langle\T I\rangle$ with respect to
$\prec_{_{T\hbox{-}gr}}$ on $\T{\B}$.\par

(iii) The set of normal monomials (mod$\langle\T I\rangle$) in
$\T{\B}$, with respect to $\prec_{_{T\hbox{-}gr}}$, is given by
$$N(\langle\T I\rangle )=\{ T^ru~|~u\in N(I),~r\in\NZ\} ,$$ where
$N(I)$ is the set of normal monomial (mod $I$) in $\B$ with respect
to $\prec_{gr}$.\vskip 6pt

{\bf Proof} The equivalence (i) $\Leftrightarrow$ (ii) is a
strengthened version of ([LWZ], Theorem 2.3.2 (i) $\Leftrightarrow$
(ii)), of which a detailed  proof was given in [LS].}\par  Noticing 
that $N(I)=\{ u\in\B~|~\LM (g){\not |} u,~g\in\G\}$, $\LM (\G)=\LM 
(\T{\G})$ where $\LM (\G )$ is taken with respect to $\prec_{gr}$ on 
$\B$ and $\LM (\T{\G})$ is taken  with respect to 
$\prec_{_{T\hbox{-}gr}}$ on $\T{\B}$, and that $\LM 
(X_iT-TX_i)=X_iT$ with respect to $\prec_{_{T\hbox{-}gr}}$ on 
$\T{\B}$, the verification of the equivalence (ii) $\Leftrightarrow$ 
(iii) is straightforward by referring to the well-known 
characterization of a Gr\"obner basis in terms of the remainder on 
division by $\T{\G}$.\QED\v5

With the preparation made above, we are ready to mention and prove
the next{\parindent=0pt\v5

{\bf 3.5. Theorem} With notation as fixed above,  let $\G =\{
g_1,\ldots ,g_s\}$ be a finite (but not necessarily homogeneous)
Gr\"obner basis for the ideal $I=\langle\G\rangle$ with respect to
$\prec_{gr}$ on $\B$, such that $\LM (\G )\cap X=\emptyset$ and $\LM
(\G )$ is reduced (in the sense of Section 2). Consider the algebra
$A=\KX /I$ which has the $\NZ$-filtration $FA$ as constructed
before.  If the monomial algebra $\OV A=\KX /\langle\LM (\G
)\rangle$ is semi-prime, then the Rees algebra $\T A$ of $A$ is
semiprimitive.\vskip 6pt

{\bf Proof} First note  that if the $\NZ$-graded lexicographic
ordering $\prec_{gr}$ on $\B$ is defined subject to
$$X_{i_1}\prec_{gr}X_{i_2}\prec_{gr}\cdots\prec_{gr}X_{i_n},$$
then the $\NZ$-graded lexicographic ordering
$\prec_{_{T\hbox{-}gr}}$ on $\T{\B}$ is defined subject to
$$T\prec_{_{T\hbox{-}gr}}X_{i_1}\prec_{_{T\hbox{-}gr}}X_{i_2}\prec_{_{T\hbox{-}gr}}\cdots\prec_{_{T\hbox{-}gr}}X_{i_n}.$$
Moreover, we also bear in mind that $d(T)=1$, and that each $X_i$
has the same positive degree as fixed in $\KX$.}\par

Let $\langle\T I\rangle$ be the non-central homogenization ideal of
$I$ in $\KXT$ with respect to $T$. Then, by Proposition 3.4(iii),
the set of normal monomials (mod$\langle\T I\rangle$) in $\T{\B}$
with respect to $\prec_{_{T\hbox{-}gr}}$ is given by $N(\langle\T
I\rangle )=\{ T^ru~|~u\in N(I),~r\in\NZ\} ,$ where $N(I)$ is the set
of normal monomial (mod $I$) in $\B$ with respect to $\prec_{gr}$.
If $T^ru_1$, $T^su_2\in N(\langle\T I\rangle )$ with
$d(T^ru_1)=d(T^su_2)$ such that
$T^ru_1\prec_{_{T\hbox{-}gr}}T^su_2$, then it follows from the
definition of $\prec_{_{T\hbox{-}gr}}$ that
$$r\ge s~\hbox{and }~u_1\prec_{gr}u_2.\eqno{(1)}$$\par

Consider the $\NZ$-gradation of $\KXT$ determined by the assigned
degrees for $T$ and $X_i$'s. If $H\in\KXT_p$ is a homogeneous
element of degree $p$, and if $H\not\in \langle\T I\rangle$, then
since $\T{\G}=\{\T g~|~g\in\G\}\cup\{ X_iT-TX_i~|~1\le i\le n\}$ is
a homogeneous Gr\"obner basis for $\langle\T I\rangle$ with respect
to $\prec_{_{T\hbox{-}gr}}$ (Proposition 3.4(ii)), the division by
$\T{\G}$ yields a representation $H=D+\sum_i\lambda_iT^{r_i}u_i$,
where $D\in\langle\T I\rangle$, $\lambda_i\in K-\{ 0\}$, and
$T^{r_i}u_i\in N(\T I)$ such that $d(T^{r_i}u_i)=p=d(H)$ for all
$i$. Put $H'=\sum_i\lambda_iT^{r_i}u_i$ and consider the
$\NZ$-graded algebra $\KXT /\langle\T I\rangle
=\oplus_{p\in\NZ}(\KXT_p+\langle\T I\rangle )/\langle\T I\rangle$.
Then the nonzero homogeneous element $\OV H$ of degree $p$
represented by $H$ in $(\KXT_p +\langle\T I\rangle)/\langle\T
I\rangle$ has the representation
$$\OV H=\OV{H'}=\sum_i\lambda_i\OV{T^{r_i}u_i},\quad \lambda_i\in K-\{ 0\},~u_i\in N(I)~
\hbox{with}~d(T^{r_i}u_i)=p=d(H).\eqno{(2)}$$
Moreover, By the above (1), if $\LM (H')=T^{r_1}u_1$ with respect to
$\prec_{_{T\hbox{-}gr}}$, then
$$\LM \left(\sum_i\lambda_iu_i\right )=u_1~\hbox{with respect to}~\prec_{gr}.\eqno{(3)}$$\par

By the definition of $\T I$, it is not difficult to verify that the
map $$\begin{array}{ccccc} \psi :&\KXT /\langle\T
I\rangle&\mapright{}{}&\KX /I=A&\\
&\OV{T}&\mapsto&1&\\
&\OV{X_i}&\mapsto&\OV{X_i}&1\le i\le n\end{array}$$ is well defined, 
and that $\psi$ is a $K$-algebra epimorphism. If,  as described in 
(2) above, $\OV H=\OV{H'}=\sum_i\lambda_i\OV{T^{r_i}u_i}$ is a 
nonzero homogeneous element of degree $p$ in $\KXT /\langle\T 
I\rangle$, then since $u_i\in N(I)$ and conventionally the ideal $I$ 
considered is a proper ideal, we have
$$\psi (\OV H)=\psi (\OV{H'})=\sum_i\lambda_i\OV{u_i}=\OV{\sum_i\lambda_iu_i}\ne 0.\eqno{(4)}$$\par

By ([LWZ], [Li1], [Li2]), $\T A\cong \KXT /\langle\T I\rangle$ as
$\NZ$-graded algebras, that is, the algebra isomorphism gives rise
to isomorphisms of $K$-vector spaces $$\T A_p\cong (\KXT_p+\langle\T
I\rangle )/\langle\T I\rangle, \quad p\in\NZ.$$ Identifying $\T A$ 
with $\KXT /\langle\T I\rangle$, we now proceed to deal with the 
Jacobson radical $J(\T A)$ of $\T A$. If $J(\T A)\ne \{ 0\}$, then  
$J(\T A)$ is a graded ideal of $\T A$ by Theorem 3.1.  Taking a 
nonzero homogeneous element of $J(\T A)$, say $\OV H\in J(\T A)\cap 
\T A_p$, where $\T A_p=(\KXT_p +\langle\T I\rangle)/\langle\T 
I\rangle$ and $H\in\KXT_p$, we may replace $\OV H$ by $\OV{H'}$ as 
in (2) above. Without loss of generality we assume that $\LM 
(H')=T^{r_1}u_1$. Our aim below is to show that 
{\parindent=.7truecm\vskip 6pt

\item{$(*)$} the normal ~monomial $u_1$ is noncyclic, thereby~
$\OV{u_1}\in J(\OV A)$ by Proposition 2.2, ~where\par $\OV{u_1}$ is
the residue class represented by $u_1$ in $\OV A$. }\vskip 6pt

Assume the contrary that $u_1$ is cyclic (see Section 2 for the
definition). Then, by Proposition 2.2, there is a monomial $v\in\B$
such that $(u_1v)^q\not\in\langle\LM (\G )\rangle$ for all 
$q\in\NZ$, or equivalently,
$$(u_1v)^q\in N(I)~\hbox{for all}~q\in\NZ ,~\hbox{in particular}~u_1v\in
N(I).\eqno{(5)}$$ On the other hand, $\OV{H'v}=\OV{H'}\OV v\in J(\T 
A)$. As $\OV{H'v}$ is again a homogeneous element of $J(A)$, it 
follows from Theorem 3.1 that
$$\OV{H'v}^m=0~\hbox{for some integer}~m>0,~\hbox{i.e.,}~(H'v)^m\in\langle\T I\rangle .\eqno{(6)}$$
Furthermore, put $H''=\sum_i\lambda_iu_iv$. By the foregoing (3),
$\LM (H'')=u_1v$ with respect to $\prec_{gr}$. Since $\psi 
(\OV{H'})\ne 0$ and $u_1v\in N(I)$ by (4) and (5), we have $\psi 
(\OV{H'v})=\OV{H''}\ne 0$ in $A$. But  it follows from (6) that 
$\OV{H''}$ is a nilpotent element in $A$, i.e.,
$$\OV{H''}^m=\left (\sum_i\lambda_i\OV{u_iv}\right )^m=\left (\psi
(\OV{H'v})\right )^m=\psi\left (\OV{H'v}^m\right )=0.\eqno{(7)}$$ 
Hence, $(H'')^m\in I$ and this gives rise to
$$(u_1v)^m=(\LM (H''))^m=\LM ((H'')^m)\in\langle\LM (I)\rangle=\langle\LM (\G )\rangle .\eqno{(8)}$$
Clearly, (8) contradicts (5). Therefore, $u_1$ is noncyclic and 
consequently $\OV{u_1}\in J(\OV A)$, proving the claim $(*)$.
\par

Finally, suppose that the monomial algebra $\OV A=\KX /\langle\LM
(\G )\rangle$ is semi-prime. Then it follows from  Theorem 2.1(ii)
that $\OV A$ is semiprimitive and hence $J(\OV {A})=\{ 0\}$. By the
above argument, we conclude that $J(\T A)=\{ 0\}$; otherwise, by the
claim $(*)$, there would be a normal monomial $u_i\in N(I)$ such 
that $\OV{u_1}\in J(\OV A)=\{ 0\}$ and hence, $u_1\in \langle\LM (\G 
)\rangle$, which is a contradiction. This shows that $\T A$ is 
semiprimitive. \QED\v5

We end this paper by the following{\parindent=0pt\v5

{\bf Open question} Let the $K$-algebra $A=\KX /\langle\G\rangle$ be
as in Theorem 3.5. If the monomial algebra $\OV A=\KX /\langle \LM
(\G )\rangle$ is semi-prime, is $A$ semiprimitive?}\v5

\centerline{References}{\parindent=1.25truecm\par
\def\hang{\hangindent\parindent}
\def\textindent#1{\indent\llap{#1\enspace}\ignorespaces}
\def\re{\par\hang\textindent}

\re{[Ber3]} G. Bergman, The diamond lemma for ring theory, {\it Adv.
Math}., 29(1978), 178--218.

\re{[G-I]} T. Gateva-Ivanova, Algorithmic determination of the
Jacobson radical of monomial algebras, in: {\it Proc. EUROCAL'85},
LNCS Vol. 378, Springer-Verlag, 1989, 355--364.

\re{[G-IL]} T. Gateva-Ivanova and V. Latyshev, On recognizable
properties of associative algebras,  {\it  J. Symbolic Computation},
6(1988), 371--388.

\re{[Gr]} E. Green, An introduction to noncommutative Gr\"obner
bases, in: {\it Computational Algebra}, Proceedings of  the fifth
meeting of the Mid-Atlantic Algebra Conference, 1993, (K. G. 
Fischer, P. Loustaunau, J. Shapiro, E. L. Green, and D. Farkas 
eds.), Lecture Notes in Pure and Applied Mathematics, Vol. 151, 
Marcel Dekker, 1994, 167--190.

\re{[Li1]} H. Li, {\it Noncommutative Gr\"obner Bases and
Filtered-Graded Transfer}, Lecture Notes in Mathematics, Vol. 1795, 
Springer, 2002.

\re{[Li2]} H. Li, $\Gamma$-leading homogeneous algebras and 
Gr\"obner bases, in: {\it Recent Developments in Algebra and Related 
Areas} (F. Li and C. Dong eds.), Advanced Lectures in Mathematics, 
Vol. 8, International Press \& Higher Education Press, 
Boston-Beijing, 2009, 155 -- 200. arXiv:math.RA/0609583, 
http://arXiv.org

\re{[Li3]} H. Li, On the calculation of gl.dim$G^{\mathbb{N}}(A)$ 
and gl.dim$\T A$ by using Gr\"obner bases, {\it Algebra Colloquium}, 
16(2)(2009), 181--194. arXiv:math.RA/0805.0686, http://arXiv.org

\re{[Li4]} H. Li, {\it Gr\"obner Bases in Ring Theory}, Monograph,
World Scientific Publishing Co., Oct. 2011.

\re{[LS]} H. Li and C. Su, On (de)homogenized Gr\"obner bases, {\it
Journal of Algebra, Number Theory: Advances and Applications},
3(1)(2010), 35--70.\\  arXiv:0907.0526, http://arXiv.org

\re{[LWZ]} H. Li, Y. Wu and J. Zhang, Two applications of
noncommutative Gr\"obner bases, {\it Ann. Univ. Ferrara - Sez. VII -
Sc. Mat.}, XLV(1999), 1--24.

\re{[Mor]} T. Mora, An introduction to commutative and
noncommutative Gr\"obner Bases, {\it Theoretic Computer Science},
134(1994), 131--173.

\re{[Row]} L.H. Rowen, {\it Ring Theory}, Vol. I, Pure and Applied
Mathematics vol. 127, Academic Press, 1988.

\re{[Uf1]} V. Ufnarovski, On the use of graphs for computing a
basis, growth and Hilbert series of associative algebras, (in
Russian 1989), {\it Math. USSR Sbornik}, 180(11)(1989), 417-428.

\re{[Uf2]} V. Ufnarovski, Introduction to noncommutative Gr\"obner
basis theory, in: {\it Gr\"obner Bases and Applications} (Linz,
1998), London Math. Soc. Lecture Notes Ser., 251, Cambridge Univ.
Press, Cambridge, 1998, 259--280.

\end{document}